\documentclass[12pt]{amsart}

\usepackage{hyperref}
\usepackage{amsthm}
\usepackage{amsmath}
\usepackage{latexsym}
\usepackage{graphicx}
\usepackage[thinlines]{easymat}
\usepackage{enumitem}
\usepackage{etoolbox}

\newcommand{\lb}{\label}
\newcommand{\f}{\frac}

\makeatletter
\patchcmd{\@settitle}{\uppercasenonmath\@title}{}{}{}
\patchcmd{\@setauthors}{\MakeUppercase}{}{}{}
\makeatother

\usepackage[letterpaper, portrait, margin=1in]{geometry}
\title{\Large\textnormal{Algebra: The Eighth Liberal Art?}}
\author{\large S. Blake Allan}
\thanks{To appear in FORMA, Winter 2025}

\begin{document}

\begin{abstract}
What is the role of algebra in classical mathematics education?  How does it relate to the four quadrivial arts?  These questions have troubled the mathematical community since the introduction of algebra into the Renaissance academy by men like François Viète, Guillame Gosselin, and René Descartes.  Their challenge is perhaps most starkly articulated at the conclusion of Viète’s Introduction to the Analytic Art, where he claims that his algebra ``appropriates to itself by right the proud problem of problems, which is: [sic] TO LEAVE NO PROBLEM UNSOLVED".
\par
Some contemporary educators respond by eschewing these methods to avoid the excessive formalization often accompanying algebra, and to give a central place to the geometrical tradition of Euclid’s Elements.  Others embrace the rise of algebra in the curriculum, focusing on contemporary techniques and priorities.
\par
This paper seeks to reconcile these perspectives by clarifying the way in which algebra participates in the quadrivial arts.  Based on testimony from both the origins of algebra and its contemporary practitioners, I argue that algebra is not so much an eighth liberal art as an arithmetical language of form – an actualized potential in arithmetic.  I conclude by offering curricular recommendations which provide glimpses of the practical insights available from this vantage.
\end{abstract}

\maketitle

\section{The Quadrivium and Vi\`ete's Challenge}	\lb{s1}

Can contemporary mathematics be reconciled with the classical {\it quadrivium}?  Whether we examine linear K-12 school curricula\footnote{Reflected in, for example, the Common Core Mathematics standards at \href{corestandards.org}{\tt corestandards.org}.} in the U.S., which try to impose a hierarchy upon loose clusters of topics, or the research triumvirate of analysis, abstract algebra,\footnote{This type of algebraic discourse is distinguished by its emphasis on formal structures – see \S\ref{s2.4} below.} and topology attempting to govern a splintering\footnote{Evidenced by the ever-expanding Mathematical Subject Classification system, the latest version of which is published openly at \href{msc2020.org}{\tt msc2020.org}.} professional discourse, the need for a robust disciplinary foundation presents itself.  It is not possible to retreat to a world lacking the profusion of present-day developments and perspectives – indeed many worthy insights would go unnoticed from such a posture.  Neither should we untether ourselves from the tradition which made many of these advances possible in the first place, as it is precisely the ends embedded in this tradition which rightly direct mathematical inquiry.  Contemporary educational practice identifies mathematics as a neutral tool to enable technological prowess, or as a fruitless exercise in symbol-pushing to distinguish high-ability students.  The belief that ``it has become necessary to change mathematics from a system of meaningful propositions into a game of formulas which is played according to certain rules"\footnote{Hermann Weyl, ``The Unity of Knowledge" [1954] in Raymond G. Ayoub, ed., {\it Musings of the Masters} (Washington, DC: The Mathematical Association of America, 2004), 74.}  has called forth a practice which is dehumanizing by virtue of its utilitarian embedded telos.  The task set for classical educators by these circumstances is thus: to reckon with the contemporary mathematical landscape in light of the commitments\footnote{One helpfully concise articulation of these commitments can be found in Brian Williams, ``Introducing Principia and Classical Education", {\it Principia: A Journal of Classical Education}, {\bf 1}, no.~1 (2022), 1-14.} shaping our distinctive educational enterprise.
\par
Recall that in Boethius’ language,\footnote{His original description in the Proemium of the {\it De Institutione Arithmetica} [c.~500 A.D.] is well worth meditating on.  A translation can be found in Michael Masi, {\it Boethian Number Theory} (New York: Rodopi, 1983), 72.} the quadrivium consists of arithmetic (concerning number in its absolute context), music (concerning harmonies among numbers), geometry (concerning fixed magnitudes), and astronomy (concerning movable magnitudes).  While these arts were neither original\footnote{Indeed, Boethius understands his project to be a Christian reception of Nicomachus of Gerasa, {\it Introduction to Arithmetic} [c.~80 A.D.], trans. Martin Luther D’ooge (London: MacMillan, 1926), and comparing their initial chapters is quite illuminating.} nor unique\footnote{One notable (and often overlooked) formulation of these arts preceding Boethius can be found in Saint Augustine, {\it On Order} [387 A.D.], trans. Michael P. Foley (New Haven: Yale University Press, 2021), 82-96.} to Boethius, his organization of them under the heading {\it quadrivium} became definitive in the medieval West.  They also extend into contemporary practice, and observing this continuity aids in clarifying the character of these arts.  Arithmetic survives under its own name in the school curriculum, but also in modern\footnote{Here and elsewhere, I understand the term ``modern" to denote a shift in Western culture beginning in the Renaissance.  Thus modern things often overlap with, but are not identical to present-day things, which are denoted ``contemporary" for clarity.} research as number theory.\footnote{Carl Friedrich Gauss’ 1801 {\it Disquistiones Arithmeticae} is commonly taken as the starting point for this area, but an interesting (pre)history is traced in Andr\'e Weil’s {\it Number Theory} (Boston: Birkh\"auser, 1984).}  Geometry also remains a pillar of the contemporary curriculum, and in modern research as topology.\footnote{There is not scope here to elucidate this connection, but mathematicians will do well to remember that an early name for topology was ``the geometry of point sets".  An accessible treatment of this remarkable area is given in Richard Earl, {\it Topology: A Very Short Introduction} (Oxford: Oxford University Press, 2019).}  Music lives on in music theory, as well as in divisibility rules and their eventual growth into ideal theory.\footnote{See John Stillwell, {\it Mathematics and Its History}, Third Edition (New York: Springer, 2010), \S 21.4.}  Astronomy has all but vanished from the contemporary mathematical curriculum, owing partially to a terminological confusion with physical accounts of the stars.\footnote{Even the contemporary English-language revival of classical education, the nature of quadrivial astronomy is not universally agreed upon.  While the Aristotelian tradition, voiced perhaps most recently in Kevin Clark and Ravi Jain, {\it The Liberal Arts Tradition}, Third Edition (Camp Hill: Classical Academic Press, 2021), Part II, places astronomy among the applied sciences (or at least the ``middle sciences"), I follow Boethius in identifying astronomy as an abstraction from celestial inquiry, just as geometry abstracts from terrestrial measurement.  Key in grounding this perspective is Plato’s remark: ``Then if, by really taking part in astronomy, we’re to make the naturally intelligent part of the soul useful instead of useless, let’s study astronomy by means of problems, as we do geometry, and leave the things in the sky alone."  Plato, {\it Republic} [c.~375 B.C.], trans. G.M.A. Grube (Indianapolis: Hackett, 1992), 530b.}  Remnants of astronomy are preserved in kinematics and in the Calculus,\footnote{In mathematical discourse, ``Calculus denotes now a certain way of performing mathematical investigations and resolutions."  Charles Hutton, {\it A Philosophical and Mathematical Dictionary}, Volume I (London: F.C. and J. Rivington, 1815), 259.  Since there are many such calculi, I reserve ``the Calculus" for the techniques of differential and integral calculus pioneered by Isaac Newton, Gottfried Leibniz, and others.} especially in the idea of mechanical curves\footnote{These are curves specifically excluded from Descartes’ early forays into algebraic geometry – see, for example, Stillwell, {\it History}, \S 13.5.} which motivated today’s theory of parametrizations.\footnote{Elaborated in, e.g., Otto Toeplitz, {\it Calculus: A Genetic Approach} [1963], trans. Luise Lange (Chicago: University of Chicago Press, 2007), Ch.~4.}
\par
The quadrivial categories underwent many challenges in the European Renaissance, but perhaps none more sustained and direct as that stemming from François Vi\`ete’s 1591 {\it Introduction to the Analytic Art}.  His purpose was to create a ``science of correct discovery in mathematics",\footnote{François Vi\`ete, {\it Introduction to the Analytic Art}, trans. T. Richard Witmer (Kent, OH: The Kent State University Press, 1983), 12.} superseding both arithmetical and geometric methods.  As Jacob Klein has elaborated upon,\footnote{In {\it Greek Mathematical Thought and the Origins of Algebra} [1934], trans. Eva Brann (Cambridge, MA: The MIT Press, 1968), Ch.~11.} Viète draws on Diophantus’ arithmetical methods, together with Pappus' paradigm of analysis and synthesis, to create a discipline he sees as ``the `one science' ($\mu\iota\alpha$ $\varepsilon\pi\iota\sigma\tau\eta\mu\eta$) which gathers all mathematical knowledge".\footnote{Klein, {\it Greek Mathematical Thought}, 159.} The Introduction ends with a powerful challenge to the rest of mathematics:
\begin{quote}
{\small
Finally, the analytic art, endowed with its three forms of zetetics, poristics, and exegetics, claims for itself the greatest problem of all, which is: NULLUM NON PROBLEMA SOLVERE.\footnote{Viète, {\it Analytic Art}, 32.  (The capitals are Viète’s).  Witmer renders the final phrase as ``To solve every problem", which while accurately capturing the grandiose character of Viète’s claim, omits the curious double-negative structure in the original.  J. Winfree Smith preserves this latter aspect, but drops the infinitive in his rendering ``There is no problem which cannot be solved" in {\it Introduction to the Analytic Art} (Annapolis: St. John’s College, 1955), 35.  Smith restores the infinitive, offering ``To leave no problem unsolved" in Klein, {\it Greek Mathematical Thought}, 353, which falls better on the ear of English speakers than the overly literal ``Nothing is not a problem to solve".}
}
\end{quote}
\par
While it is important to consider this claim in its original scope, my primary interest will be its pedagogical implications, which reverberate from Viète’s time into our contemporary discourse.  This avenue was opened by Guillaume Gosselin,\footnote{Gosselin is comparatively obscure in contemporary mathematics, but details of his academic station can be found in Natalie Zemon Davis, ``Mathematicians in Sixteenth-Century French Academies: Some Further Evidence", Renaissance News, {\bf 11}, no.~1 (1958): 3-10.} who in his 1577 {\it On the Great Art} introduced algebra into the university curriculum, and ``spoke of it as the eighth liberal art".\footnote{Smith, {\it Introduction to the Analytic Art}, 3}  This is a place of highest honor for Gosselin, as he says that algebra ``by the ancients was called the science of creation and creatures, and by others the rule of rules, and finally by others the queen of the sciences".\footnote{{\it Huius scienti\ae\  qu\ae\  ab antiquis appellata est scientia creatur\ae\ \& creaturarum, ab aliis regula regularum, ab aliis denique regina scientiarum}, Gvlielmi Gosselini, {\it De Arte Magna} (Paris, Apud Aegidium Beys, via Iacobæa, ad insigne Lilij albi, 1577), Cap.~III.}  This begs the question: Can algebra be reconciled to the traditional liberal arts in our curricula, or does it subsume them in its aspirations to universality?
\par
In what follows, I will establish that no conflict need exist between these modern creations and the ancient quadrivium.  This is the case because {\it algebra is the arithmetical language of form}.  By ``arithmetical language", I mean an organized system of discrete signs\footnote{Weyl’s concept of {\it aufweisbar} comes very close to this notion – see {\it Unity of Knowledge}, 73-75.} and terms,\footnote{This definition of ``language" is perhaps too inclusive, as German, C++, and knitting shorthand could all be called ``languages" in this sense.  A more precise (but markedly less accessible) term might be ``semiotic", but these finer points belong to a different occasion.} generalized from those used in arithmetic.  By ``form", I mean ``the ordering principle or arrangement of a given thing",\footnote{Phillip J. Donnelly, {\it The Lost Seeds of Learning} (Camp Hill: Classical Academic Press, 2021), 69.} whether it be an equation, a curve, or some collection of mathematical entities.  On this account, algebra is therefore a way of addressing the form of mathematical entities and their relationships in an arithmetical mode, whether the objects themselves are discrete or continuous.

\section{A Timeline of the Arithmetical Language}	\lb{s2}

The long and varied history of algebra provides ample resources to substantiate this view of its fundamental character.  Four key periods illustrate the role algebra plays in the liberal arts: the era preceding Viète’s {\it Introduction}, the early years of algebraic geometry (also called analytic geometry), a transition period populated by interpreters, and the contemporary structural era.  For each time period, I will attend to the most significant testimony regarding the nature of algebra itself (particularly from European authors), both providing support for the role of algebra as an arithmetical language of form, and clarifying what this means in the subdiscipline(s) occupying many\footnote{An exhaustive treatment is not fitting here, but many more details may be found in Bartel L. van der Waerden, {\it A History of Algebra: From al-Khwārizmī to Emmy Noether} (New York: Springer, 1985).} investigators in these eras.

\subsection{The Pre-Vi\`ete Era}	\lb{s2.1}

While many interesting studies\footnote{See Stillwell, {\it History}, Ch.~6 and especially Bartel L. van der Waerden, {\it Geometry and Algebra in Ancient Civilizations} (Berlin: Springer, 1983), Ch.~3.} have been conducted regarding prefigurations of algebraic thought prior to 800 A.D., what is most important for the present purpose begins with Muhammed al-Khwārizmī’s {\it The Book of Restoration and Opposition}.  Though the author delineates his subject matter as ``concerning arithmetical and geometrical problems",\footnote{Louis Charles Karpinski, {\it Robert of Chester’s Latin Translation of the Algebra of Al-Khowarizmi} (New York: MacMillan, 1915), 67} his method is decidedly arithmetical.  Many of the problems involve what we now write as quadratic equations, whose solutions are established by both arithmetical and geometrical\footnote{Karpinski, {\it Algebra}, 77ff, 129ff.} means.  However, the attention to numerical characteristics of the problem at hand distinguishes even al-Khwārizmī’s geometrical figures from the unarithmetized tradition of Euclid (which was also well-known in the Arabic-speaking world at this time).
\par
The 1145 Latin edition of {\it The Book} transliterates, rather than translates the title, coining the term ``algebra" in the sense used in the West.  The Arabic term \includegraphics[scale=.025]{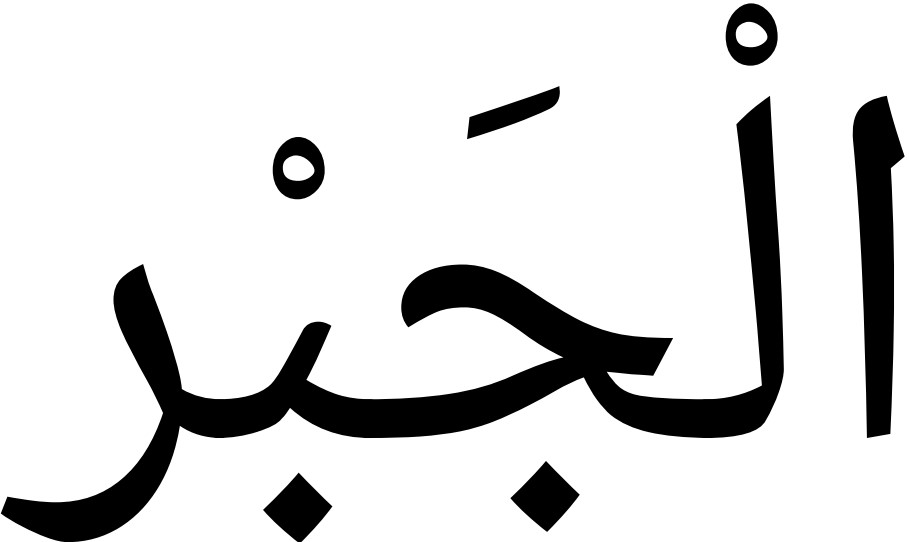} ({\it al-jabr}, restoration) indicates ``the transference of negative terms" to the opposite side of an equation, and \includegraphics[scale=.025]{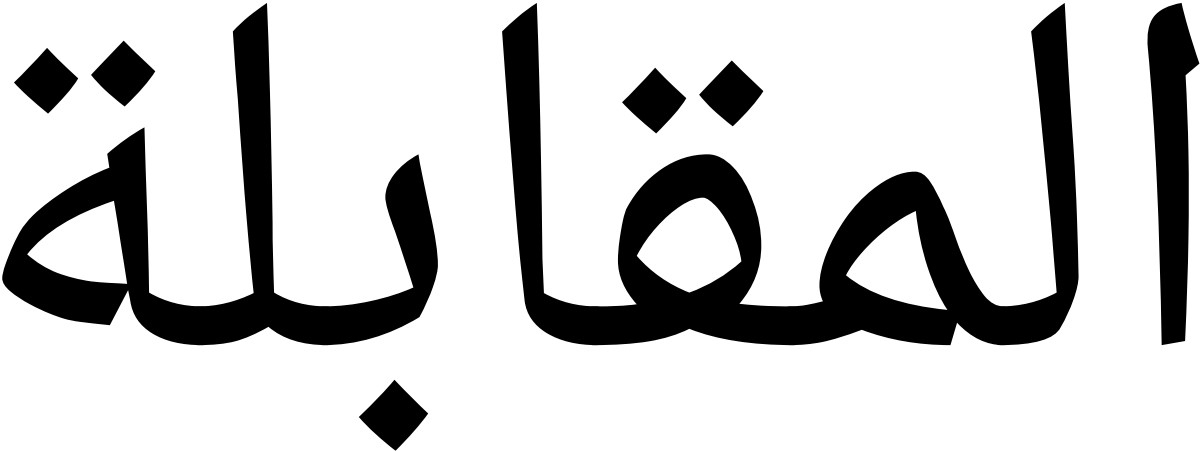} ({\it al-muqabala}, opposition) denotes ``the combination of like terms".\footnote{Karpinski, {\it Algebra}, 67.}  In al-Khwārizmī’s work these are clearly arithmetical processes which attend to the form of equations, even when the quantities involved may represent geometrical entities.
\par
Perhaps the best primary summary\footnote{Admittedly, many important contributions are left out (or are unattributed) in Cardano’s work.  Some of the gaps may be filled by, e.g., David Eugene Smith, {\it The History of Mathematics}, Volume I (New York: Dover, 1958), Ch.~VI-VII and Simon Gindikin, {\it Tales of Mathematicians and Physicists}, trans. Alan Shuchat (New York: Springer, 2007), 1-26.} of the progress of European algebra in the intervening 400 years is Girolamo Cardano’s {\it The Great Art}, published in 1545.  This immense collection of rules (general procedures) for solving equations represents substantial contributions to the theory of equations.  Particularly significant are the first prominent appearances of solutions to the cubic and quartic equations – results which do not seem to have been known before the modern period.  Cardano declares his work to be a ``most abstruse and unsurpassed treasury of the entire [subject of] arithmetic",\footnote{Girolamo Cardano, {\it The Great Art or the Rules of Algebra}, trans. T. Richard Witmer (Cambridge, MA: The MIT Press, 1968), 1.} clearly locating his work within our familiar liberal art.
\par
Indeed, the era of algebra preceding Viète’s 1591 work may be characterized as concrete, calculational, and centered on equations.  This is an enduring image of the discipline – as late as 1866, Joseph Serret attested that ``Algebra is, properly speaking, the Analysis of equations",\footnote{Quoted in Nicolas Bourbaki, {\it Elements of the History of Mathematics} [1974], trans. John Meldrum (Berlin: Springer, 1994), 52.} and in 1970, the Bourbaki committee opened their volume on the subject by declaring ``Algebra is concerned essentially with calculating".\footnote{Nicolas Bourbaki, {\it Elements of Mathematics}, Algebra I [1970] (New York: Springer, 1989), xxi.  Thomas Hobbes makes an even more general comment, remarking that ``in what matter soever there is a place for {\it addition} and {\it substraction}, there is also place for {\it Reason}; and where these have no place, there {\it Reason} has nothing at all to do." {\it Leviathan} [1651] (London: J.M. Dent \& Sons, 1914), 18.}  Reading Viète’s work in this larger context, it is not surprising that we find it to possess ``inherently numerical characteristics" with an object that ``is, its generality notwithstanding, `arithmetically' determined".\footnote{Klein, {\it Greek Mathematical Thought}, 184.}

\subsection{The Era of Algebraic Geometry}	\lb{s2.2}

One of the reasons Viète’s {\it Introduction} is commonly overlooked is the foundational role it played in enabling René Descartes' famous {\it Geometry} of 1637.  Though it has been said that this work ``reduced the methods of geometry to calculations performed on numerical quantities",\footnote{James Clerk Maxwell, {\it A Treatise on Electricity and Magnetism}, Volume I (Oxford: Clarendon Press, 1873), 9.} Descartes himself saw his notions of adding and subtracting lines as introducing ``arithmetical terms into geometry, for the sake of greater clearness."\footnote{David Eugene Smith and Marcia L. Latham, {\it The Geometry of René Descartes} (New York: Dover, 1954), 5.}  The project of his {\it Geometry} was, fittingly, geometrical, but the method decidedly arithmetical.  In particular, Descartes drew on the equational vision of algebra still dominant in his time, insisting that ``all points of those curves which we may call `geometric', that is, those which admit of precise and exact measurement, must bear a definite relation to all points of a straight line, and that this relation must be expressed by means of a single equation".\footnote{Smith and Latham, {\it Geometry}, 48.}  Today such curves are called algebraic precisely because points on them satisfy an algebraic equation whose attributes indicate the character of the curve.
\par
One example of the fruits of this method is a classification of curves by the degree of the associated equation, which is defined as the largest sum of exponents in a single term.  By reformulating ancient results, Descartes realized that straight lines were the only curves of degree $1$, and the conic sections (circle, ellipse, hyperbola, parabola) were the curves of degree $2$.\footnote{The latter family also contains some exceptional cases, such as two perpendicular lines – see Vladimir I Arnold, {\it Real Algebraic Geometry}, trans. Gerald G. Gould and David Kramer (New York: Springer, 2013), Ch.~1 \& 2.}  In 1695, Isaac Newton classified 72 species of degree-$3$ curves using novel algebraic transformations and equivalences, giving an early demonstration of the power of Descartes' approach.\footnote{And also providing an opportunity for further inquiry, as 6 species were overlooked!  Detailed references regarding these developments appear in Stillwell, {\it History}, \S 7.4.}
\par
While geometry was never absent from early algebraic developments,\footnote{For example, the complex numbers owe their heritage to a problem whose solution was geometrically apparent but algebraically subtle – see Tristan Needham, {\it Visual Complex Analysis}, 25th Anniversary Edition (Oxford: Oxford University Press, 2023), \S 1.1.} this second period emphasized the insights algebraic methods could bring to continuous magnitudes.  The persistence of this vision in the contemporary curriculum needs no elaboration, but one important echo of this perspective was offered by Hermann Grassmann in his 1862 {\it Extension Theory}, which treated extensive magnitudes by arithmetical means\footnote{See, e.g., Grassmann’s foreword in {\it Extension Theory}, trans. Lloyd C. Kannenberg (Providence: American Mathematical Society and London Mathematical Society, 2000).} and inaugurated the area now known as linear algebra.\footnote{An contemporary introduction to this field can be found in Sheldon Axler, {\it Linear Algebra Done Right}, Third Edition (New York: Springer, 2015).}

\subsection{The Era of Interpreters}	\lb{s2.3}

Just as Leonhard Euler became the inheritor of Gottfried Leibniz’s calculus, so too did he receive and transform the algebraic tradition.  Many of his early works elaborate what later became known as number theory,\footnote{C. Edward Sandifier, {\it The Early Mathematics of Leonhard Euler} (Providence: MAA Press, 2007) is full of interesting examples in this vein.} and his vision of the Calculus centered around an especially algebraic version of the function concept.\footnote{See Gindikin, {\it Tales}, 199-203.}  Most pertinent to the present discussion is his {\it Elements of Algebra}, published in 1770.  In the first chapter, Euler says that in algebra, ``we consider only numbers, which represent quantities, without regarding the different kinds of quantity."\footnote{Leonhard Euler, {\it Elements of Algebra}, Fifth Edition, trans. John Hewlett (London: Longman, 1840), 2.}  At first glance this seems to limit algebra only to numerals, but he immediately clarifies that arithmetic ``treats of numbers in particular, and is the {\it science of numbers properly so called}".\footnote{Euler, {\it Algebra}, 2, emph. original.}  This shows (together with the content of the {\it Elements of Algebra}, which does not limit itself only to calculation regarding particular numbers) that Euler sees algebra as an arithmetical way to investigate any quantity appearing in mathematics.  In fact, he contends that ``the foundation of all the Mathematical Sciences must be laid in a complete treatise on the science of numbers, and in an accurate examination of the different possible methods of calculation.  This fundamental part of mathematics is called Analysis, or Algebra."\footnote{Euler, {\it Algebra}, 2.}  This parallel to Boethius’ account of arithmetic as the mathematical art ``which holds the principal place and position of mother to the rest"\footnote{Masi, {\it Boethian Number Theory}, 74} clearly locates algebra as participating in the arithmetical tradition.  As hinted by the title’s reference to Euclid, Euler’s {\it Elements} lives up to this expansive vision, weaving together elementary number theory, the solutions of equations (quadratic, cubic, and quartic), roots, logarithms, proportions, infinite series, and much more.
\par
The most significant extension of Euler's explicit principles was given in George Peacock's 1830 {\it Treatise on Algebra} (later expanded into two volumes).  Besides recounting many important results formulated by Euler and his successors, Peacock articulates a very general division within algebra, which helps clarify the relationship between Viète's innovations and the more general algebra of Viète's inheritors.  Work focused on equations and their solutions (even in general form) falls under Peacock's category of arithmetical algebra, ``the science which results from the use of symbols and signs to denote numbers and the operations to which they may be subjected\ldots being used in the same sense and with the same limitations as in common arithmetic".\footnote{George Peacock, {\it A Treatise on Algebra, Volume I: Arithmetical Algebra} (Cambridge: Cambridge University Press, 1842), 1.}  On the other hand, the more general inquiries of Descartes, Euler, and others belong to symbolical algebra.  Peacock grounds this latter discipline on what he calls the principle of the permanence of equivalent forms, which states that ``whatever algebraical forms are equivalent, when the symbols are general in form but specific in value, will be equivalent likewise when the symbols are general in value as well as in form."\footnote{George Peacock, {\it A Treatise on Algebra, Volume II: Symbolical Algebra} (Cambridge: Cambridge University Press, 1845), 59.}  As an example of this principle, consider the general multiplicative property $x^{m}\cdot x^{n}=x^{m+n}$ for powers. This is certainly true when $m$ and $n$ are positive integers, but Peacock's principle asserts its truth when $m$ and $n$ are any two numbers.  By choosing, say, $m,\,n=\f{1}{2}$, we learn that $x^{1/2}\cdot x^{1/2}=x$, a property already known to hold for the square root of $x$.  Thus, we can discover the familiar equivalence $\sqrt{x}=x^{1/2}$ between fractional exponents and radicals.
\par
Euler and Peacock’s\footnote{John Dubbey makes an interesting argument in {\it The Mathematical Work of Charles Babbage} (Cambridge: Cambridge University Press, 1978), Ch.~5 that Peacock’s insights in this direction may have been partially or wholly anticipated by Charles Babbage in the latter’s unpublished manuscript ``The Philosophy of Analysis".  However, due to the tragic cyber-attack on the British Library in October 2023, I have been unable to examine Babbage’s work in full to evaluate Dubbey’s claim.} articulations of algebraic principles, and their thorough use of these principles in their work provide an important foundation for an account of algebra as the arithmetical language of form.  As a closing demonstration of the value of these principles, William Rowan Hamilton’s change of opinion comes to mind.  Hamilton initially held to a curious conviction that algebra was the ``science of pure time", and attempted to reduce the complex numbers to pairs of real numbers using this principle.\footnote{The results of this strange attempt are elaborated in William Rowan Hamilton, ``Theory of Conjugate Functions, or Algebraic Couples" [1837], in {\it The Mathematical Papers of Sir William Rowan Hamilton}, Volume III (Cambridge: Cambridge University Press, 1967), 3-96, and critiqued by Arthur Cayley, {\it Presidential Address to the British Association, September 1883}, in William Ewald, {\it From Kant to Hilbert}, Volume I (Oxford: Clarendon Press, 1996), Ch.~14.}  In later work, Hamilton admits the peculiarity of his own approach,\footnote{In an impressively long footnote in ``Preface to `Lectures on Quaternions'" [1853] in {\it Mathematical Papers Volume III}, 125-126.} and through the lens of Peacock’s views obtains a much clearer organization of the various imaginary systems (complex numbers, quaternions, etc.) under his consideration. This instance demonstrates that not only are our interpreters' views philosophically coherent, they also lead to substantial mathematical insights.

\subsection{The Era of Structures}	\lb{s2.4}

The final epoch in algebra includes much of contemporary discourse, but its foundations are interlaced with the previous era's insights.  This structural era was touched off by a peculiar discovery: Niels Abel's corrected attempt to solve the quintic equation in radicals.  In 1824, he announced that such a general solution, which had been found rapidly in the sixteenth century for the cubic and quartic, was not possible at all!\footnote{This statement is sometimes abbreviated to ``the unsolvability of the quintic", but key question to ask is: By what means?  A solution ``in radicals" is a finite formula for the root(s) of an equation using only addition, subtraction, multiplication, division, and extraction of roots (such as the quadratic formula).  While no such formula exists for an arbitrary quintic equation, there is an explicit solution in terms of more complicated objects – see J. V. Armitage and W. F. Eberlein, {\it Elliptic Functions} (Cambridge: Cambridge University Press), Ch.~10.}  Évariste Galois brought Abel’s result\footnote{Some details of which can be found in Jeremy Gray, {\it A History of Abstract Algebra} (New York: Springer, 2018), \S 9.3.} into a systematic theory of the (un)solvability of algebraic equations in 1831, but Galois’ sudden death\footnote{The story of Galois’ brief life and profound contributions can be found in van der Waerden, {\it A History}, Ch.~6.} delayed his work’s reception into the broader mathematical community.  His investigation hinged on permutations of the roots of an equation, which when studied as a collection became the paradigmatic example of a mathematical group.  Examples of groups include rotations, translations, reflections, and permutations.  Each is equipped with a multiplication or composition that joins two members of the group, a concept which ``seems inseparable from the first rudiments of calculation"\footnote{Bourbaki, {\it Elements of History}, 47.} in Bourbaki’s view.
\par
Equation theory motivated the group concept,\footnote{However, the theory of groups can also be compellingly motivated through the theory of Euclidean symmetries – see Needham, {\it VCA}, \S 1.4.} and considerations of number and divisibility gave rise to other more complicated algebraic structures.  Detailed descriptions of these more exotic structures – fields, rings, and many more – exceed our present scope,\footnote{Stillwell, {\it History}, Ch.~19-21 offers an accessible discussion.} but the clear emphasis of this era is on {\it collections} of objects, rather than the individuals belonging to them.  Groups have become the standard expression of the symmetries of an object or arrangement,\footnote{As narrated in Hermann Weyl, {\it Symmetry} (Princeton: Princeton University Press, 1952).} while ideals within rings offer factorizations of prime numbers in strange domains.\footnote{This is discussed, for example, in David Hilbert, {\it The Theory of Algebraic Number Fields} [1897], trans. Iain T. Adamson (New York: Springer, 1998), Ch.~2.}  Over time, similarities among these varied structures led Samuel Eilenberg and Saunders Mac Lane to formulate category theory,\footnote{Mac Lane's own discussion of this idea can be found in his volume {\it Mathematics: Form and Function} (New York: Springer, 1986), \S XI.9.} a conceptual detachment of algebraic form from the particular structures that embody it.  This method continues to inform contemporary treatments of algebra, and lead Mac Lane to describe the overall discipline by saying ``Algebra tends to the study of the explicit structure of postulationally defined systems closed with respect to one or more rational operations."\footnote{Mac Lane, “Some Recent Advances in Algebra”, The American Mathematical Monthly, {\bf 46}, no.~1 (1939): 3-19}
\par
In light of algebra's position as the arithmetical language of form, it is thus the formal character that receives the most emphasis in the contemporary era.  The arithmetical character is evident in the notation and some of the methodology, but these are seen as secondary by many contemporary practitioners.  In any case, the structural approach to algebra taken by the last century and a half of innovators continues to underscore the depth of our view.

\section{Pedagogical Reflections}	\lb{s3}

In the previous section, we have seen that our definition of algebra is both motivated by and accounts for a wide range of historical testimony about the discipline - what implications might this have for contemporary mathematics curricula?
\par
One aid brought by these considerations is clarity regarding the place of algebra in relation to curricular goals.  While the linear order of contemporary school mathematics often places algebra as the final stepping stone in the race to the Calculus, recognizing its participation in the arithmetical tradition helps open broader possibilities for algebra courses.  Several such possibilities have already been realized in historical usage of Euler's {\it Elements} – for instance, Section I of Part I would form a helpful capstone to any arithmetic curriculum, while Section III would serve well in a course on proportion and music.  Selections from Sections II and IV would provide more than enough for a typical algebra course, with plenty of supplementary material for especially curious students.
\par
Furthermore, once it is clear that algebra is not all-encompassing, the scope of what it {\it does} usefully illuminate and the character of this illumination becomes much clearer.  The application of algebra to geometry has been of prime importance since the reception of Descartes, but as Ravi Jain has recently argued,\footnote{In {\it The Enchanted Cosmos} (Camp Hill, PA: Classical Academic Press, 2025), Ch.~3.} it is also essential to maintain a robust category of the discrete.  In contrast to standard coordinate techniques, Josiah Williard Gibbs' vector arithmetic\footnote{Accessibly exposited in the first two chapters of Edwin Bidwell Wilson, {\it Vector Analysis} (Cambridge, MA: Yale University Press, 1901).} and Felix Klein’s account of symmetry\footnote{Beautifully illustrated for a general audience in David Mumford, Caroline Series, and David Wright, {\it Indra's Pearls: The Vision of Felix Klein} (Cambridge: Cambridge University Press, 2002).} maintain and utilize this distinction to obtain remarkable insights.  Students of logic frequently benefit from George Boole’s algebraic system which represents a proposition ``by an equation the form of which determines the rules of conversion and transformation, to which the given proposition is subject."\footnote{George Boole, {\it The Mathematical Analysis of Logic} (Cambridge: MacMillan, Barclay, \& MacMillan, 1847), 8.}  Finally, we mention that even so deep a result as Abel's unsolvability theorem of 1824 can be both understood and established by upper school students,\footnote{Vladimir Arnold gave a course with this aim 1963-64, and the resulting notes were published as V.B. Alekseev, {\it Abel’s Theorem in Problems and Solutions} [1976], trans. Francesca Aicardi (Dordrecht: Kluwer Academic Publishers, 2004).} drawing on an algebraic encapsulation of ideas from Johannes Kepler’s study of the icosahedron in his 1619 {\it The Harmony of the World}.

\newpage
\section{Conclusion}	\lb{s4}

As we have seen, testimony from the entire history of algebra in the West supports its identification as the arithmetical language of form – an actualized potential in the traditional liberal art of arithmetic.  This expansive definition of the discipline accounts for both its computational and structural aspects, identifying it as an important (but not all-encompassing) item for inclusion in classical curricula.  It is my hope that this short survey can provide a better acquaintance with this valuable art, and another avenue for attending to the One who ``hast ordered all things in measure, and number, and weight".\footnote{Wis.~11:20 (KJV).}

\end{document}